\documentclass[a4paper,11pt, reqno]{amsart}
\usepackage{array}
\usepackage{graphicx}
\usepackage{float}
\usepackage{amsmath,latexsym}
\usepackage{amssymb}
\usepackage{geometry}
\usepackage{amsfonts}
\usepackage{hyperref}
\usepackage{color}
\usepackage{multirow}

\newcommand{\dsum}{\displaystyle\sum}

\usepackage{setspace}
\allowdisplaybreaks

\let\origmaketitle\maketitle
\def\maketitle{
  \begingroup
  \def\uppercasenonmath##1{} 
  \let\MakeUppercase\relax 
  \origmaketitle
  \endgroup
}

\def\R{\mathbb{R}}
\def\Z{\mathbb{Z}}

\usepackage{soul}

\begin{document}

\title[]{\huge Multiclass Optimal Classification Trees with SVM-splits}

\author[V. Blanco, A. Jap\'on \MakeLowercase{and} J. Puerto]{{\large V\'ictor Blanco$^\dagger$, Alberto Jap\'on$^\ddagger$ and  Justo Puerto$^\ddagger$}\medskip\\
$^\dagger$Institute of Mathematics (IMAG), Universidad de Granada\\
$^\ddagger$Institute of Mathematics (IMUS), Universidad de Sevilla\\
\texttt{vblanco@ugr.es}, \texttt{ajapon1@us.es}, \texttt{puerto@us.es}
}

\date{\today}

\maketitle 

\begin{abstract}
In this paper we present a novel mathematical optimization-based methodology to construct tree-shaped classification rules for multiclass instances. Our approach consists of building Classification Trees in which, except for the leaf nodes, the labels are temporarily left out and grouped into two classes by means of a SVM separating hyperplane. We provide a Mixed Integer Non Linear Programming formulation for the problem and report the results of an extended battery of computational experiments to assess the performance of our proposal with respect to other benchmarking classification methods.
\end{abstract}

\section{Introduction}

Interpretability is a crucial requisite demanded to machine learning methods provoked by the tremendous amount of methodologies that have arised in the last decade~\cite{du2019techniques}. It is expected that the model that results when applying a machine learning methodology using a training sample, apart from being able to adequately predict the behaviour of out-of-sample observations, can be interpreted. Different tools have been applied to derive \textit{interpretable} machine learning methods. One of the most popular strategies to simplify the obtained models is \textit{feature selection}, in which a reduced set of attributes is to be selected without loosing quality in the predictions. Reducing the number of parameters to analyze, the models can be easier to understand, yielding higher descriptive accuracy. One could also consider models that can be modulated, in the sense that a great proportion of its prediction-making process can be interpreted independently. This is the case of generalized linear models \cite{hastie2017generalized}. Other methods incorporate interpretability as a synonym of being able to be reproduced by humans in its entire construction~\cite{ct,letham2015interpretable}. This is the case of Decision Trees with small depth which can be visualized and interpreted easily by users even not familiar with the tools behind their construction. We adopt a tree-based methodology through this paper.

Among the wide variety of strands derived under the lens of Machine Learning, classification is one that has attracted a lot of attention because of its applicability in many different fields~\cite{writing,credit,insurance,cancer,cleveland}. Classification methodologies aim to adequately predict the class of new observations provided that a given sample has been used to construct the \textit{classification rule}. The role of Mathematical Programming in the construction of classification models has been widely recognized, and some of the most popular methods to derive classification rules are based on solving optimization problems \cite{bpr20,svm,OCTsurvey,BJP20,gunluk2018optimal}. Moreover, Mathematical Programming has also been proven to be a flexible and accurate tool when requiring interpretability to the obtained models ~\cite{baldomero2020tightening,baldomero2021robust,BJP21,gaudioso2017lagrangian}. 

However, most of the optimization tools derived to construct classifiers assume instances with only two classes. In this paper, we provide a novel classification method in which the instances are allowed to be classified into two or more classes. The method is constructed using one of the most interpretable classification method, Classification Trees, but combined with Support Vector Machines, which provide highly predictive models.

We have developed a Mathematical Programming model that allows to construct an Optimal Classification Tree for a given training sample, in which each split is generated by means of a SVM-based hyperplane. When building the tree, the labels of the observations are ignored in the branch nodes, and they are only accounted for in the leaf nodes where misclassification errors are considered. The classification tree is constructed  to minimize the complexity of the tree (assuring interpretability) and also the misclassification risk (assuring predictive power).

\subsection{Related Works}

Several machine learning methodologies have been proposed in the literature in order to construct highly predictive classification rules. The most popular ones are based on Deep Learning mechanisms~\cite{Agarwal_2018}, $k$-Nearest Neighborhoods \cite{knn1,knn}, Na\"ive Bayes \cite{nb}, Classification Trees (CT)  \cite{ct,friedman2001} and Support Vector Machines (SVM)~\cite{svm}.  Among them, CT and SVM, which are, by nature, optimization-based methodologies, apart from producing highly predictive classifiers have been proven to be very flexible tools since both allow the incorporation of different elements (through the adequate optimization models by means of constraints and objective functions) to be adapted to different situations, as Feature Selection~\cite{gunluk2018optimal,baldomero2020tightening,baldomero2021robust,jimenez2021novel}, accuracy requirements~\cite{benitez2019cost,gan2021robust} or dealing with unbalanced or noisy instances~\cite{eitrich2006efficient,BJP20_SVM,blanquero2021optimal}, amongst others.

Support Vector Machines were originally introduced by Cortes and Vapnik~\cite{svm} as a binary classification tool that builds the decision rule by means of a separating hyperplane with large separation between the two classes. This hyperplane is obtained by solving a convex quadratic optimization problem, in which the goal is to separate data by their two differentiated classes, maximizing the margin between them and minimizing the misclassification errors. Duality properties of this optimization problem allow one to extend the methodology to find nonlinear separators by means of kernels. Classification  Trees were firstly introduced by Breiman et. al \cite{ct}, and the decision rule is based on a hierarchical relation among a set of nodes which is used to define paths that lead observations from the root node (highest node in the hierarchical relation), to some of the leaves in which a class is assigned to the data. These paths are obtained according to different optimization criteria over the predictor variables of the training sample. The decision rule comes up naturally, the classes predicted for new observations are the ones assigned to the terminal nodes in which observations fall in. Clearly, the classification rules derived from CTs are easily interpretable by means of the splits that are constructed at the tree nodes. In \cite{ct}, a greedy heuristic procedure, the so-called CART approach, is presented to construct CTs. Each level of the tree is sequentially constructed: starting at the root node and using the whole training sample, the method minimizes an impurity measure function obtaining as a result a split that divides the sample into two disjoint sets which determine the two descendant nodes. This process is repeated until a given termination criteria is reached (minimum number of observations belonging to a leaf, maximum depth of the tree, or minimum percentage of  observations of the same class on a leaf, amongst others). In this approach, the tree  grows following a top-down greedy approach, an idea that is also shared in other popular decision tree methods like C4.5~\cite{quinlan93} or ID3~\cite{quinlan96}. The advantage of these methods is that the decision rule can be obtained rather quickly even for large training samples, since the whole process relies on  solving manageable problems at each node. Nevertheless, these types of heuristic approaches may not obtain the \textit{optimal} classification tree, since they look for the best split locally at each node, not taking into account the splits that will come afterwards. Thus, these local branches may not capture the proper structure of the data, leading to misclassification errors in out-of-sample observations. Furthermore, the solutions provided by these methods can result into very deep (complex) trees, resulting in overfitting and, at times, loosing interpretability of the classification rule. This difficulty is usually overcome by pruning the tree as it is being constructed by comparing the gain on the impurity measure reduction with  respect to the complexity cost of the tree. 

The recent advances on modeling and solving difficult Optimization problems  together with the flexibility and adaptability of these models have motivated the use of optimization tools to construct supervised classification methods with a great success~\cite{bertsimas2019machine,OCTsurvey}). In particular, recently, Bertsimas and Dunn \cite{bertsimas2017optimal} introduced the notion of \textit{Optimal Classification Trees} (OCT) by approaching Classification and Regression Trees under optimization lens, providing a Mixed Integer Linear Programming formulation for its optimal construction.  Moreover, the authors proved that this model can be solved for reasonable size datasets, and equally important, that for many different real datasets, significant improvements in accuracy with respect to CART can be obtained. In contrast to the standard CART approach, OCT builds the tree by solving a single optimization problem taking into account (in the objective function) the complexity of the tree,  avoiding post pruning processes. Moreover, every split is directly applied in order to minimize the misclassification errors on the terminal nodes, and hence, OCTs are more likely to capture the hidden patterns of the data.  

While SVM were initially designed to deal only with bi-class instances, some extensions have been proposed in the literature for multiclass classification. The most popular multiclass SVM-based approaches are One-Versus-All (OVA) and One-Versus-One (OVO). The former, namely OVA, computes, for each class $r \in \{1, \ldots, k\}$, a binary SVM classifier labeling the observations as $1$, if the observation is in the class $r$, and $-1$ otherwise. The process is repeated for all classes ($k$ times), and then each observation is classified into the class whose constructed hyperplane is the furthest from it in the positive halfspace. In the OVO approach, classes are separated with ${k}\choose{2}$ hyperplanes using one hyperplane for each pair of classes, and the decision rule comes from a voting strategy in which the most represented class among votes becomes the class predicted. OVA and OVO inherit most of the good properties of binary SVM. In spite of that,  they are not able to correctly classify datasets where separated clouds of observations may belong to the same class (and thus are given the same label) when a linear kernel is used.  Another popular method is the directed acyclic graph SVM, DAGSVM~\cite{Agarwal_2018}. In this technique,  although the decision rule involves the same hyperplanes built with the OVO approach, it is not given by a unique voting strategy but for a sequential number of votings in which the most unlikely class is removed until only one class remains. In addition, apart from OVA and OVO, there are some other methods based on decomposing the original multiclass problem into several binary classification ones. In particular,  in \cite{Allwein} and \cite{Bakiri}, this decomposition is based on the construction of a coding matrix that determines the pairs of classes that will be used to build the separating hyperplanes. Alternatively, other methods such as CS~(\cite{cs}), WW~(\cite{ww}) or LLW~(\cite{llw}), do not address the classification problem sequentially but as a whole considering all the classes within the same optimization model. Obviously, this seems to be the correct approach. In particular, in WW, $k$ hyperplanes are used to separate the $k$ classes, each hyperplane separating one class from the others, using $k-1$ misclassification errors for each observation. The same separating idea, is applied in CS but reducing the number of misclassification errors for each observation to a unique value. In LLW, a different error measure is proposed to cast the Bayes classification rule into the SVM problem implying theoretical statistical properties in the obtained classifier. These properties cannot be ensured in WW or CS.

We can also find a quadratic extension based on LLW proposed by \cite{Guermeur}. In \cite{genSVM}, the authors propose a multiclass SVM-based approach,  \textit{GenSVM}, in which the classification boundaries for a problem with $k$ classes are obtained in a $(k-1)$-dimensional space using a simplex encoding.
Some of these methods have become popular and are implemented in  most software packages in machine learning as \texttt{e1071} \cite{e1071}, \texttt{scikit-learn} \cite{python} or \cite{MSVMpack}. Finally, in the recent work \cite{BJP20} the authors propose an alternative approach to handle multiclass classification extending the paradigm of binary SVM classifiers by construnting a polyhedral partition of the feature space and an assignment of classes to the \emph{cells} of the partition, by maximizing the separation between classes and minimizing two intuitive misclassification errors. 

\subsection{Contributions}

In this paper, we propose a novel approach to construct Classification Trees for multiclass instances by means of a mathematical programming model. Our method is based on two main ingredients: (1) An optimal binary classification tree (with oblique cuts) is constructed in the sense of \cite{bertsimas2017optimal}, in which the splits and pruned nodes are determined in terms of the misclassification errors at the leaf nodes; (2) The splits generating the branches of the tree are build by means of  binary SVM-based hyperplanes separating fictitious clases (which are also decided by the model), i.e., maximizing separation between classes and minimizing the distance-based misclassification errors. \\
Our specific contributions include:
\begin{enumerate}
\item Deriving an interpretable multiclass classification rule which combines two of the most powerful tools in supervised classification, namely OCT and SVM.
\item The classifier is constructed using a mathematical programming model that can be formulated as a Mixed Integer Second Order Cone Programming problem. The classifier is simple to apply and interpretable.
\item Several valid inequalities are presented for the formulation that allow one to strengthen the model and to solve larger size instances in smaller CPU times.
\item An extensive battery of computational experiments on realistic datases from UCI is reported showing that our approach outperforms other decision tree-based methodologies as CART, OCT and OCT-H.
\end{enumerate}

\subsection{Paper structure}

Section \ref{sec:2} is devoted to fix the notation and to recall the tools that are used to derive our method. In Section \ref{sec:3}
 we detail the main ingredients of our approach and illustrate its performance on a toy example. The mathematical programming model that allows us to construct the classifier is given in Section \ref{sec:4}, where we include all the elements involved in the model: parameters, variables, objective function and constraints. In Section \ref{sec:5} we report the results of our experiments to assess the performance of our method compared with other tree-shaped classifiers. Finally, Section \ref{sec:6} is devoted to draw some conclusions and future research lines on the topic.
 
\section{Preliminaries}\label{sec:2}

This section is devoted to introduce the problem under study and to fix the notation used through this paper. We also recall the main tools involved in our proposed approach namely, Support Vector Machines and Optimal Classification Trees. These methods are adequately combined to develop a new method, called Multiclass Optimal Classification Trees with Support Vector Machines based splits (MOCTSVM).

We are given a training sample, $\mathcal{X} = \left\lbrace (x_1,y_1),\ldots ,(x_n,y_n)\right\rbrace \subseteq  \mathbb{R}^p \times \left\lbrace 1,\ldots , K \right\rbrace$, which comes up as the result of measuring $p$ features over a set of $n$ observations $(x_1,\ldots , x_n)$ as well as a label in $\{1, \ldots, K\}$ for each of them $(y_1,\ldots , y_n)$. The goal of a classification method is to build a decision rule so as to accurately assign labels ($y$) to data ($x$) based on the behaviour of the given training sample $\mathcal{X}$.

The first ingredient that we use in our approach is the Support Vector Machine method. SVM is one of the most popular optimization-based methods to design a classification rule in which only two classes are involved, usually referred as the positive $(y=+1)$ and the negative class $(y=-1$). The goal of linear SVM is to construct a hyperplane separating the two   classes by maximizing their separation and simultaneously minimizing the misclassification and margin violation errors.  Linear SVM can be formulated as the following convex optimization problem:
\begin{align}
\min & \ \frac{1}{2}\|\omega\|_2^2 + c\dsum_{i\in N} e_i &\nonumber\\
\mbox{s.t.} & \;\; y_i(\omega'x_i+\omega_0) \geq 1-e_i,   & \forall i\in N, \nonumber \\
&  \omega \in \mathbb{R}^p, \ \omega_0 \in \mathbb{R},  & \nonumber\\
& e_i \in \mathbb{R}_+, & \forall i\in N. \nonumber
\end{align}
where $c$ is the regularization parameter that states the trade-off between training errors and model complexity (margin), $\omega'$ is the transpose of the vector $\omega$ and $\|\cdot\|_2$ is the Euclidean norm in $\R^p$ (other norms can also be considered but still keeping similar structural properties of the optimization problem~\cite{bpr20}). Note that with this approach,  the positive (resp. negative) class will tend to lie on the positive (resp. negative) half space induced by the hyperplane $\mathcal{H} = \{z \in \R^p: \omega' z+\omega_0 =0\}$. On the other hand, the popularity of SVM is mostly due to the so call kernel trick.  This allows one to project the data onto a higher dimensional space in which a linear separation is performed in a most accurate way with no need of knowing such a space, but just knowing the form of its inner products, and  maintaining the computational complexity of the optimization problem (see \cite{svm} for further details).


The second method that we combine in our approach is Classification Trees. CTs are a family of classification methods based on a hierarchical relationship among a set of nodes. These methods allow one to create a partition of the feature space by means of hyperplanes that are sequentialy built. CT starts on a node containing the whole sample, that is called the root node, in which the first split is applied. When applying a split on a node (by means of a hyperplane separating the observations, two new branches are created leading to two new nodes, which are referred to as its child nodes. The nodes are usually distinguished into two groups: branch nodes, that are nodes in which a split is applied, and on the other hand the leaf nodes, which are the terminal nodes of the tree. Given a branch node and a hyperplane split in such a node, their branches  (left and right) are defined as each of the two halfspaces defined by the hyperplane. The final goal of CT is to construct branches in order to obtain leaf nodes as pure as possible with respect to the classes. In this way, the classification rule for a given observation consists of assigning it to the most popular class of the leaf where it belongs to.

\begin{figure}[H]
\begin{center}
\fbox{\includegraphics[scale=0.06]{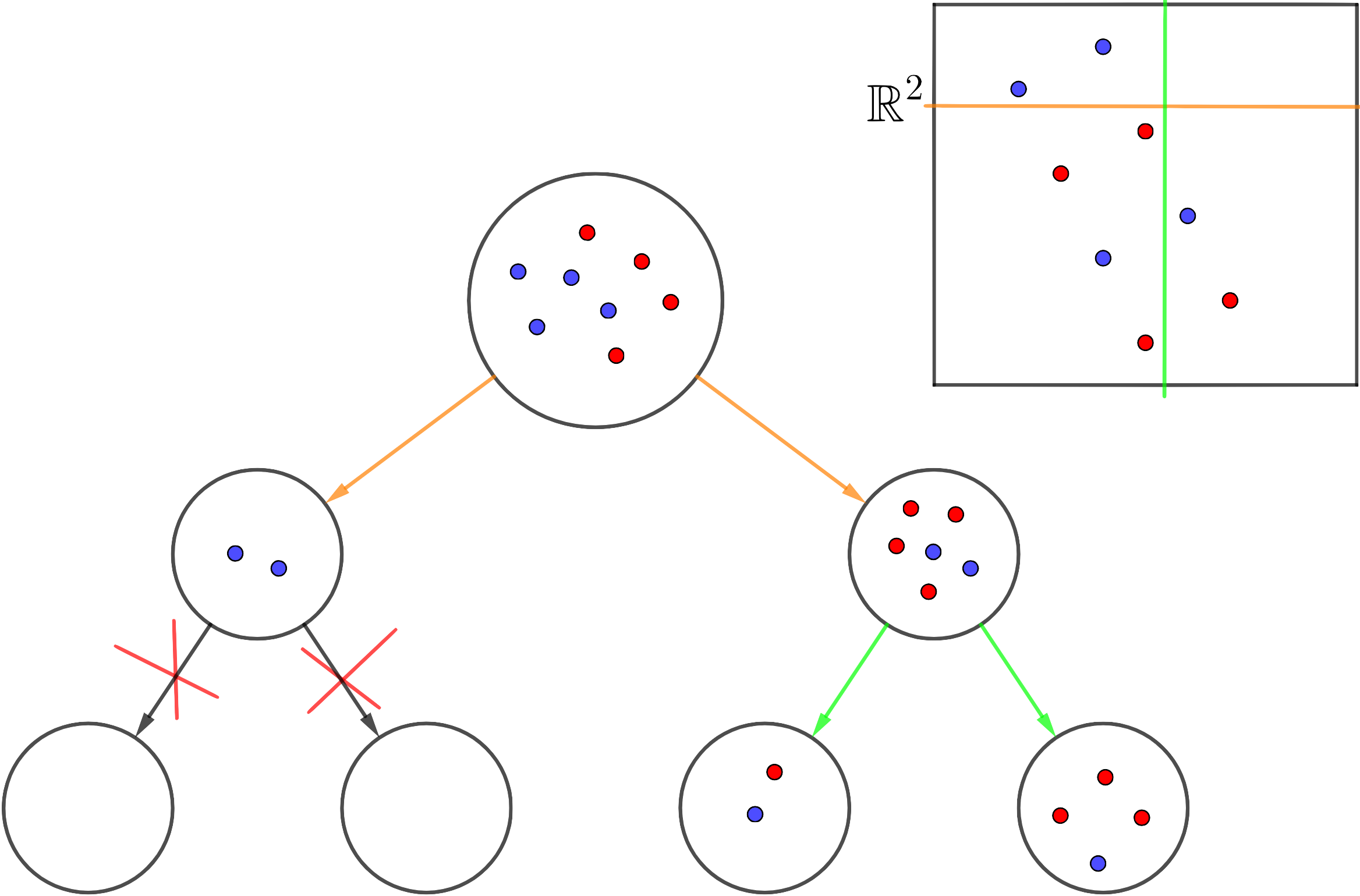}}~\fbox{\includegraphics[scale=0.06]{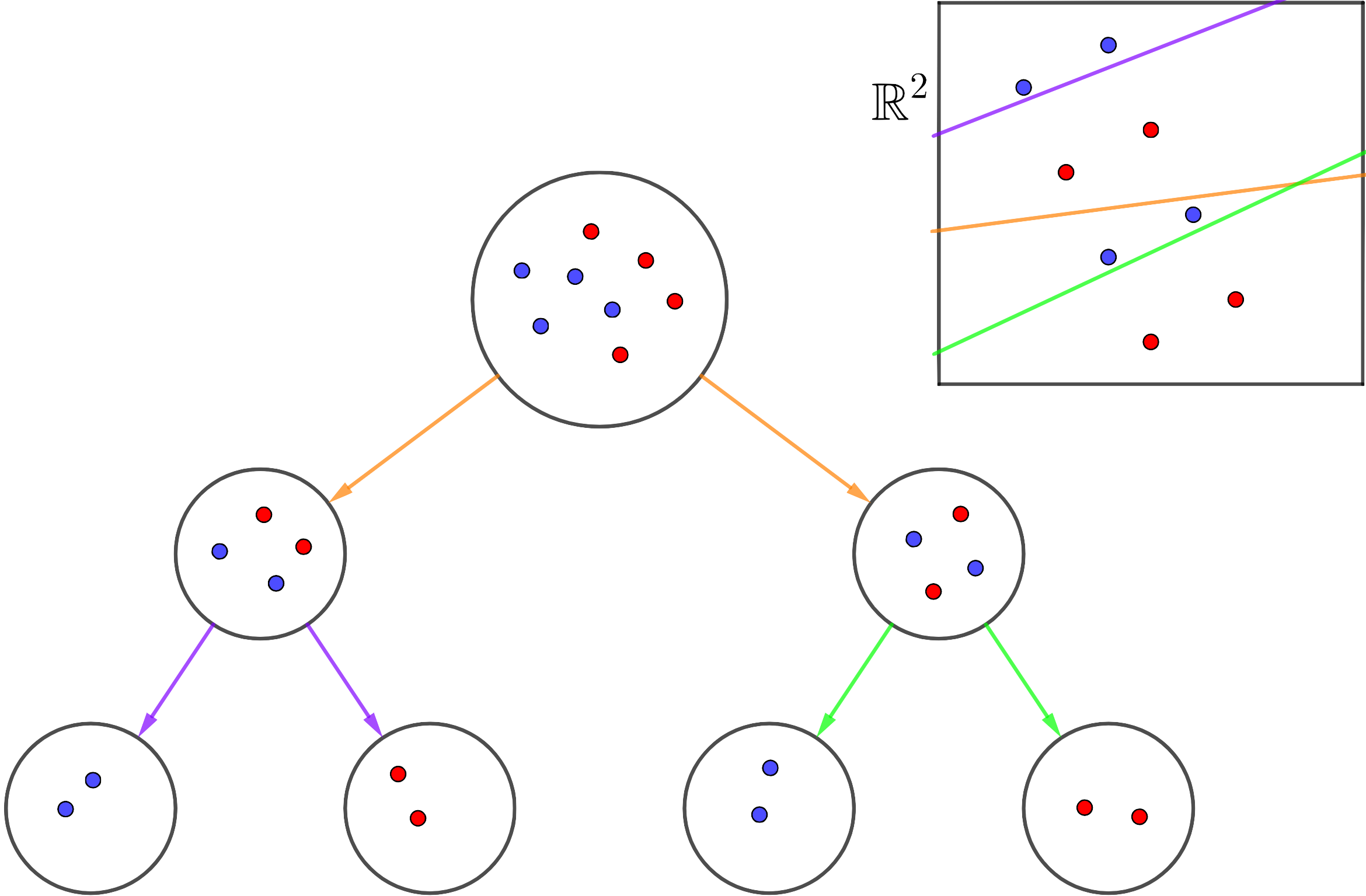}}
\caption{Example of a CT obtained with CART (left) and OCT-H (right) approaches for the same instance.}\label{fig:1}
\end{center}
\end{figure}

There is a vast amount of literature on CTs since they provide an  easy interpretable classification rule. One of the most popular methods to construct CT is known as CART, introduced in \cite{ct}. CART is a greedy heuristic approach that starts at the root node looking for the split in a single feature that minimizes a given impurity function, creating two new nodes. The same procedure is sequentially applied until a stop criteria is reached (maximal depth of the tree, proportion of observations of a single class in a node, etc). The main advantage of CART is its low computational cost, since nowadays very deep trees can be obtained within a few seconds. However, CART does not guarantee the \emph{optimality} of the classification tree, in the sense that more accurate trees could be obtained if instead of locally constructing the branches one looks at the final configuration of the leaf nodes. For instance, in Figure \ref{fig:1}(left) we show a CT constructed by CART for a biclass problem with maximal depth $2$. We draw the classification tree, and also in the top right corner, the partition of the feature space (in this case $\R^2$). As can be observed, the obtained classification is not perfect (not all leaf nodes are composed by pure classes) while in this case is not difficult to construct a CT with no classification errors. This  situation is caused by the myopic construction done by the CART approach that, at each node only cares on better classification at their children, but not at the final leaf nodes, while subsequent branching decisions clearly affect the overall shape of the tree.

Motivated by this drawback of CART, in \cite{bertsimas2017optimal}, the authors propose an approach to build an Optimal Classification Tree (OCT) by solving a single Mathematical Programming problem in which not only single-variable splits are possible but \emph{oblique} splits involving more than one predictive variable (by means of general hyperplanes in the feature space) can be constructed. In Figure \ref{fig:1}(right) we show a solution provided by OCT with hyperplanes (OCT-H) for the same example. One can observe that when splitting the root node (orange branches) a good local split is not obtained (the nodes contain half of the observations in different classes), however, when adding the other two splits, the final leaves only have observations of the same class, resulting in a perfect classification rule for the training sample.

\begin{figure}[h]
\begin{center}
\includegraphics[scale=0.08]{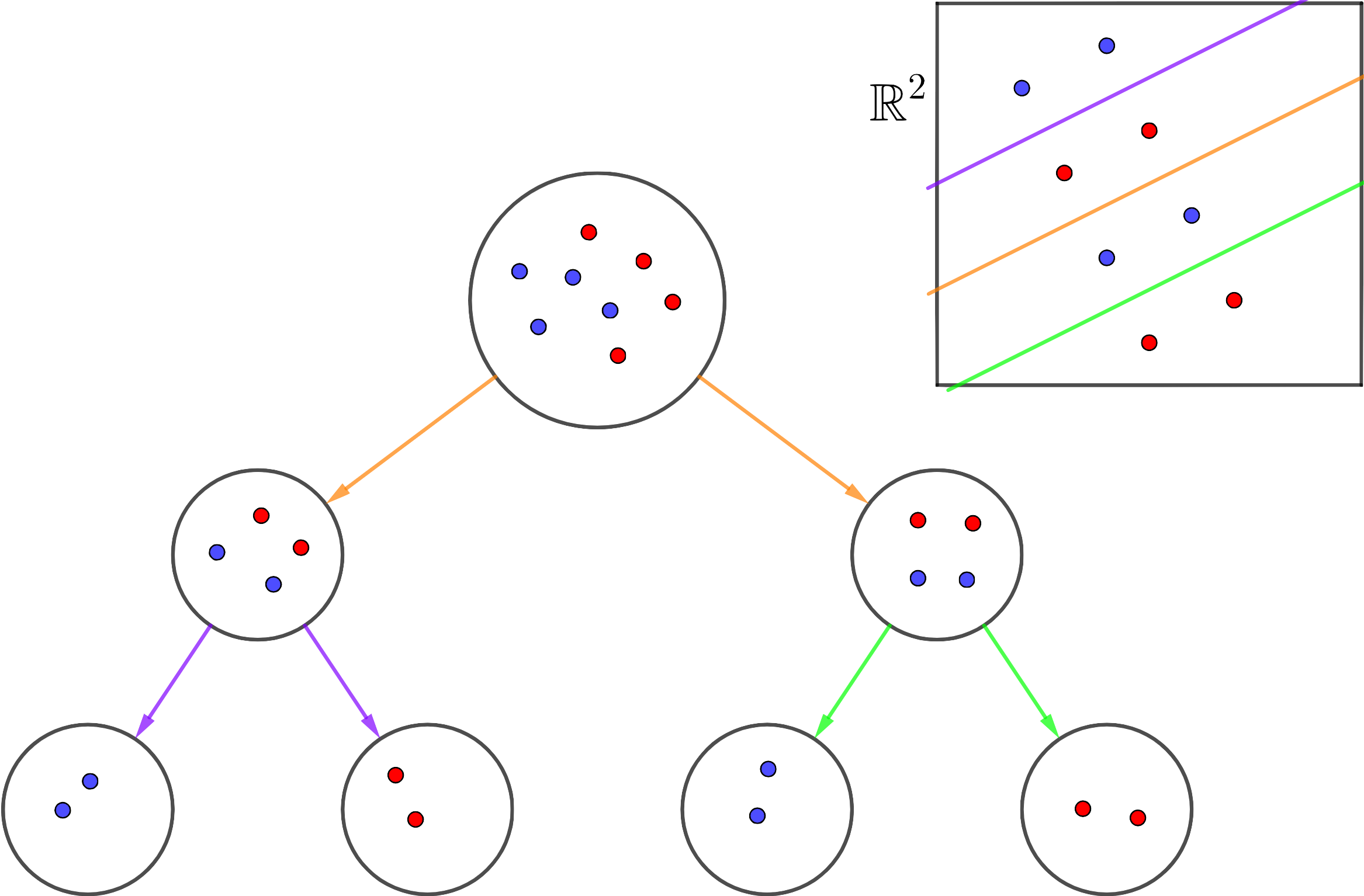}~
\caption{Example of a CT obtained with OCTSVM.}\label{fig:2}
\end{center}
\end{figure}

Both approaches, OCTs and SVMs can be combined in order to construct classification trees in which the classes separated by the hyperplanes determined in the CT are \emph{maximally} separated, in the sense of the SVM approach. This idea is not new and has been proven to outperform standard optimal decision trees methods amongst many different biclass classification problems, as for instance, in \cite{BJP21} where the OCTSVM method is proposed. In Figure \ref{fig:2} we show how one could construct OCTs with larger separations between the classes using OCTSVM but still with the same $100\%$ accuracy in the training sample as in OCT-H, but more protected to misclassification in out-sample observations.

Nevertheless, as far as we know, the combination of OCT and SVM has only been analyzed for biclass instances. The extension of this method to multiclass settings  (more than two classes) is not trivial, since one could construct more complex trees or use a multiclass SVM-based methodology (see e.g. \cite{cs,ww,llw}).  However, these adaptations of the classical SVM method have been proved to fail in real-world instances (see \cite{BJP20}).  In the rest of the paper we describe a novel methodology to construct accurate multiclass tree-shaped classifiers based on a different idea: constructing CTs with splits induced by bi-class SVM separators in which the classes of the observations at each one of the branch nodes are determined by the model, but adequately chosen to provide small classification errors at the leaf nodes. The details of the approach are given in the next section.

\section{Multiclass OCT with SVM splits}\label{sec:3}

In this section we describe the method that we propose to construct classification rules for multiclass instances, in particular Classification Trees in which splits are generated based on the SVM paradigm.

\begin{figure}[h]
\centering
\includegraphics[scale=0.75]{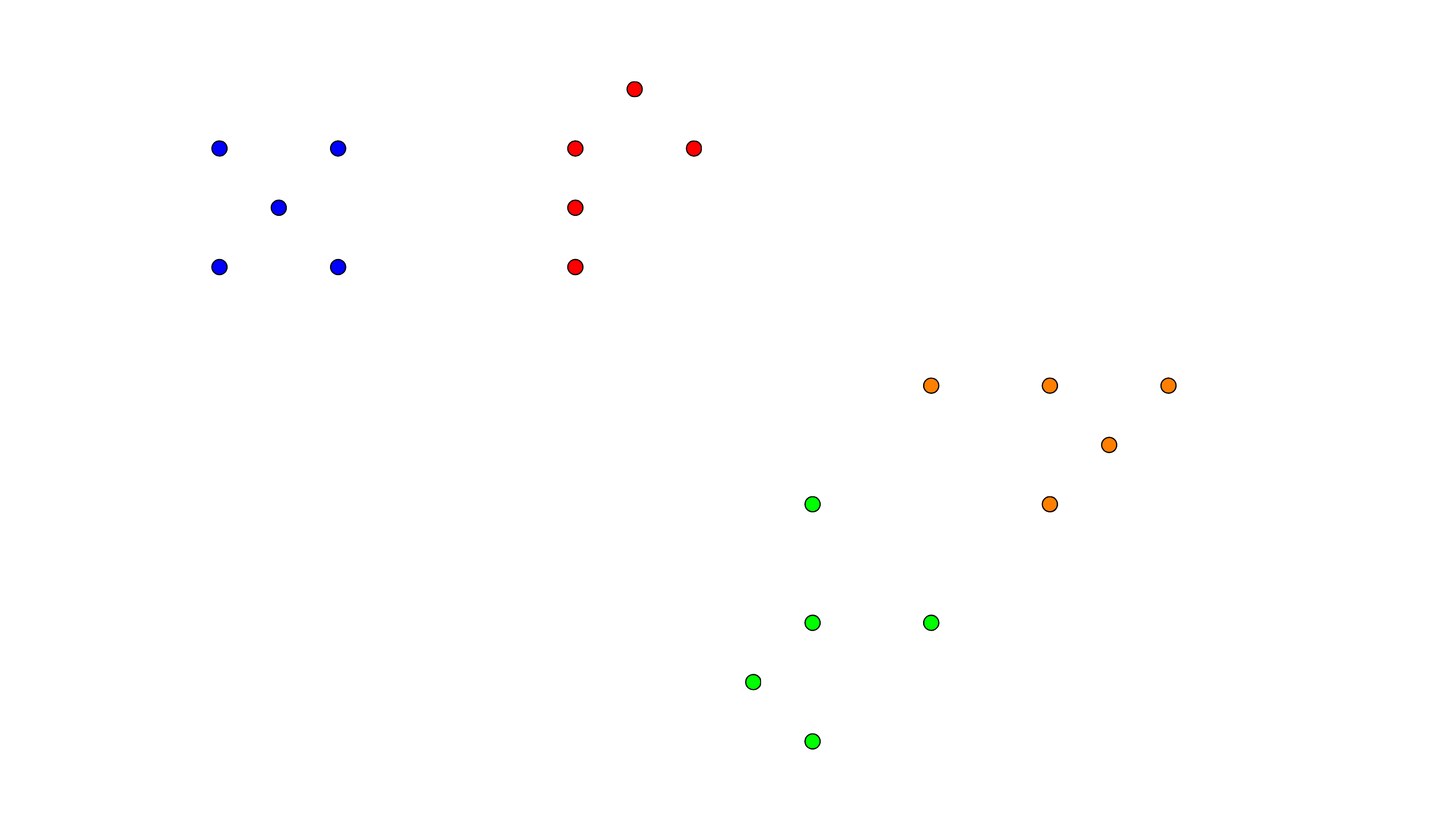}
\caption{Instance for a $4$-class problem.}
\label{fig:3}
\end{figure}

As already mentioned, our method is based on constructing OCT with SVM splits, but where the classes of the observations are momentarily ignored and only accounted for at the leaf nodes. In order to illustrate the idea under our method, in Figure \ref{fig:3} we show a toy instance with a set of points with four different classes (blue, red, orange and green).

\begin{figure}[h]
\centering
\includegraphics[scale=0.55]{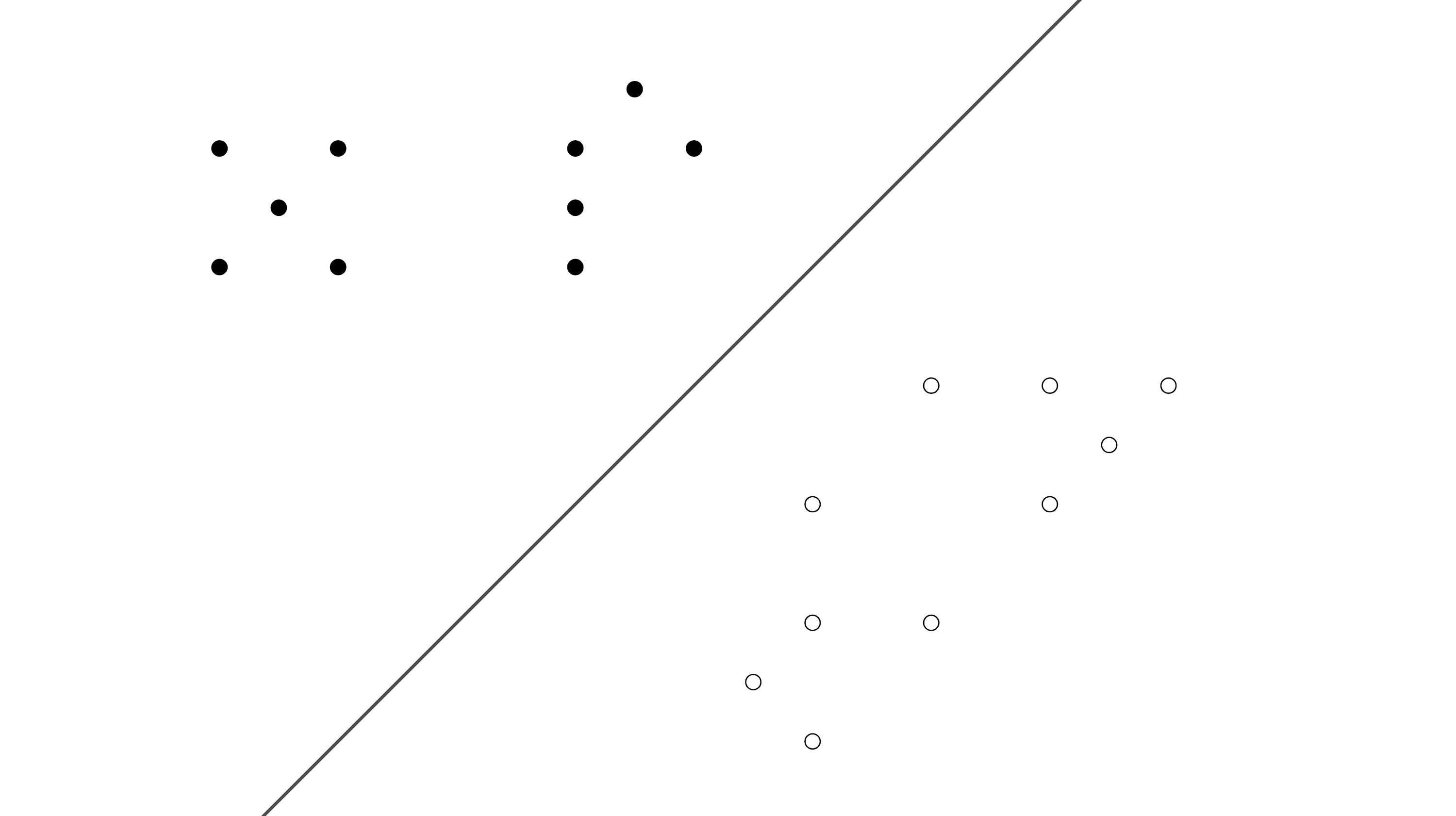}
\caption{Root split on the 4-class classification problem}\label{fig:4}
\end{figure}

\begin{figure}[h]
\centering
\includegraphics[scale=0.55]{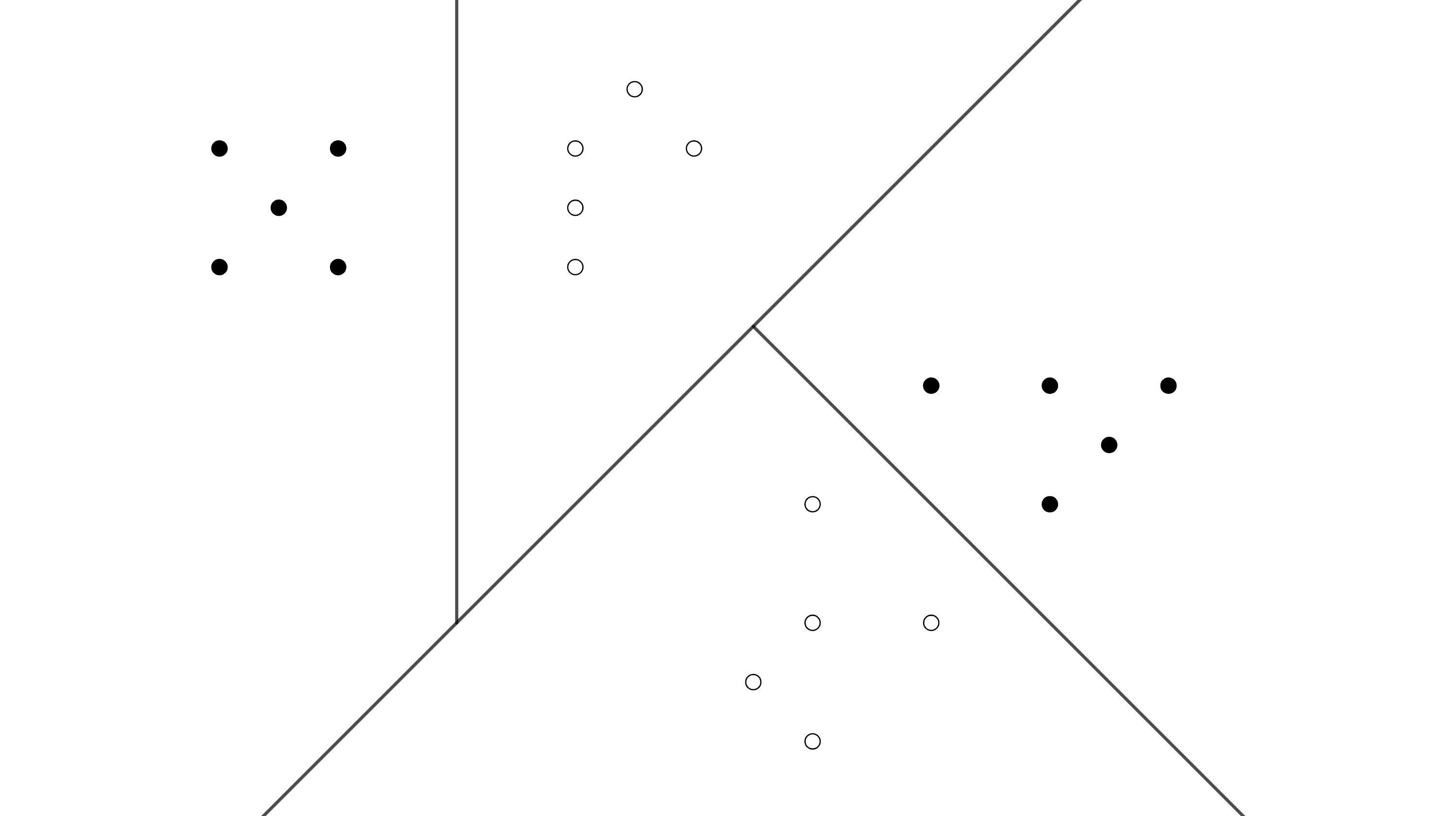}
\caption{Child node splits on the 4-class classification problem higheleted as the fictitious classes decided in our model.}\label{fig:5}
\end{figure}

First, at the root node (the one in which all the observations are involved), our method constructs a SVM separating hyperplane for two fictitious classes (which have to be also decided). A possible separation could be the one shown in Figure \ref{fig:4}, in which the training dataset has been classified into two classes (black and white). 
This separation allows one to generate two child nodes, the black and the white nodes. At each of these nodes, the same idea is applied until the leaf nodes are reached. In Figure \ref{fig:5} we show the final partition of the feature space according to this procedure. 

\begin{figure}[h]
\centering
\includegraphics[scale=0.55]{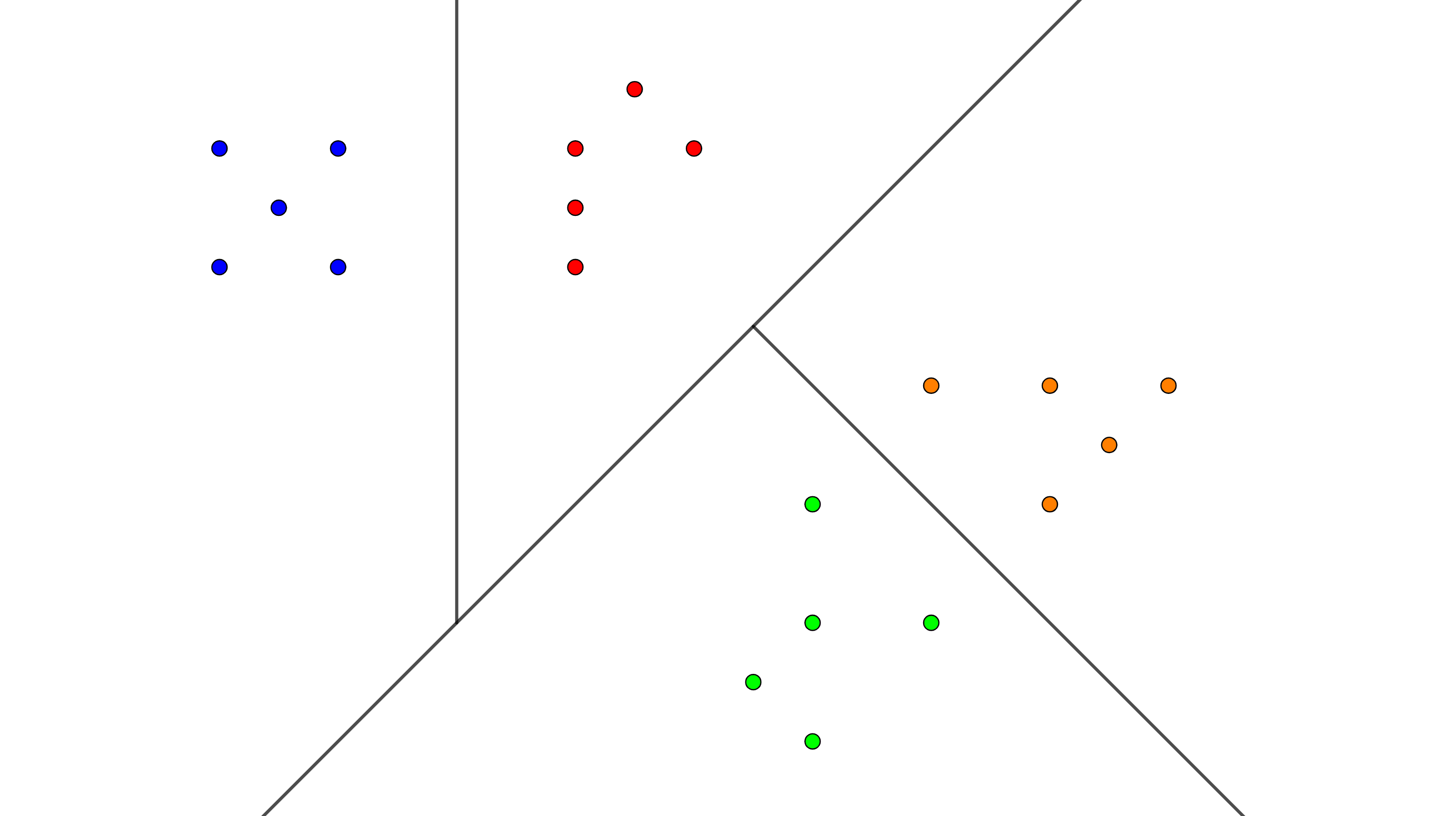}
\caption{Child node splits on the 4-class classification problem with their original labels (colors).}\label{fig:6}
\end{figure}

Clearly, ignoring the original classes of the training sample in the whole process would result in senseless trees, unless one accounts for the goodness in the classification rule in the training sample at the leaf nodes. Thus, at the final leaf nodes, the original labels are recovered and the classification is performed according to the generated hyperplanes. The final result of this tree is shown in Figure \ref{fig:6} where one can check that the constructed tree achieves a perfect classification of the training sample.

Once the tree is constructed with this strategy, the decision rule comes up naturally as it is usually done in decision trees methods, that is, out of sample observations will follow a path on the tree according to the splits and they will be assigned to the class of the leaf where they lie in (the most represented class of the leaf over the training set). In case a branch is pruned when building the tree, observations will be assigned to the most represented class of the node where the prune took place.

\section{Mathematical Programming Formulation for MOCTSVM}\label{sec:4}

In this section we derive a Mixed Integer Non Linear Programming formulation for the MOCTSVM method described in the previous section.

We assume to be given a training sample $\mathcal{X} = \left\lbrace (x_1,y_1),\ldots ,(x_n,y_n)\right\rbrace \subseteq  \mathbb{R}^p \times \left\lbrace 1,\ldots , K \right\rbrace$. We denote by $N=\{1, \ldots, n\}$ the index set for the observations in the training sample. We also consider the binary representation of the labels $y$ as:
 $$
 Y_{ik} = \left\{\begin{array}{cl}
 1 & \mbox{if $y_i=k$,}\\
 0 & \mbox{otherwise,}
 \end{array}\right. \mbox{for all}\ i\in N,\ k = 1,\ldots K .
 $$
 Moreover, without loss of generality we will assume the features to be normalized, i.e., $x_1,\ldots , x_n \in [0,1]^p$.
 
We will construct decision trees with a fixed maximum depth $D$. Thus, the classification tree is formed by at most $T=2^{D+1}-1$ nodes. We denote by $\tau=\{1, \ldots, T\}$ the index set for the tree nodes, where node $1$ is the root node and nodes $2^D, \ldots, 2^{D+1}-1$ are the leaf nodes.

For any node $t\in \tau\backslash\{1\}$, we denote by $p(t)$ its (unique) parent  node. The tree nodes can be classified in two sets: branching and leaf nodes. The branching nodes, that we denote by $\tau_b$, will be those in which the splits are applied.  In constrast, in the leaf nodes, denoted by $\tau_l$, no splits are applied but is where predictions take place. The branching nodes can be also classified into two sets: $\tau_{bl}$ and $\tau_{br}$ depending on whether they follow the left or the right branch on the path from their parent nodes, respectively. $\tau_{bl}$ nodes are indexed with even numbers meanwhile $\tau_{br}$ nodes are indexed with odd numbers.

We define a level as a set of nodes which have the same depth within the tree. The number of levels in the tree to be built is $D+1$ since the root node is assumed as the zero-level. Let $U = \{u_0, \ldots, u_D\}$ be the set of levels of the tree, where each $u_s \in U$ is the set of nodes at level $s$, for $s=0, \ldots, D$. With this notation, the root node is $u_0$ while $u_D$ represent the set of leaf nodes.

\begin{figure}[h]
\includegraphics[width=\textwidth]{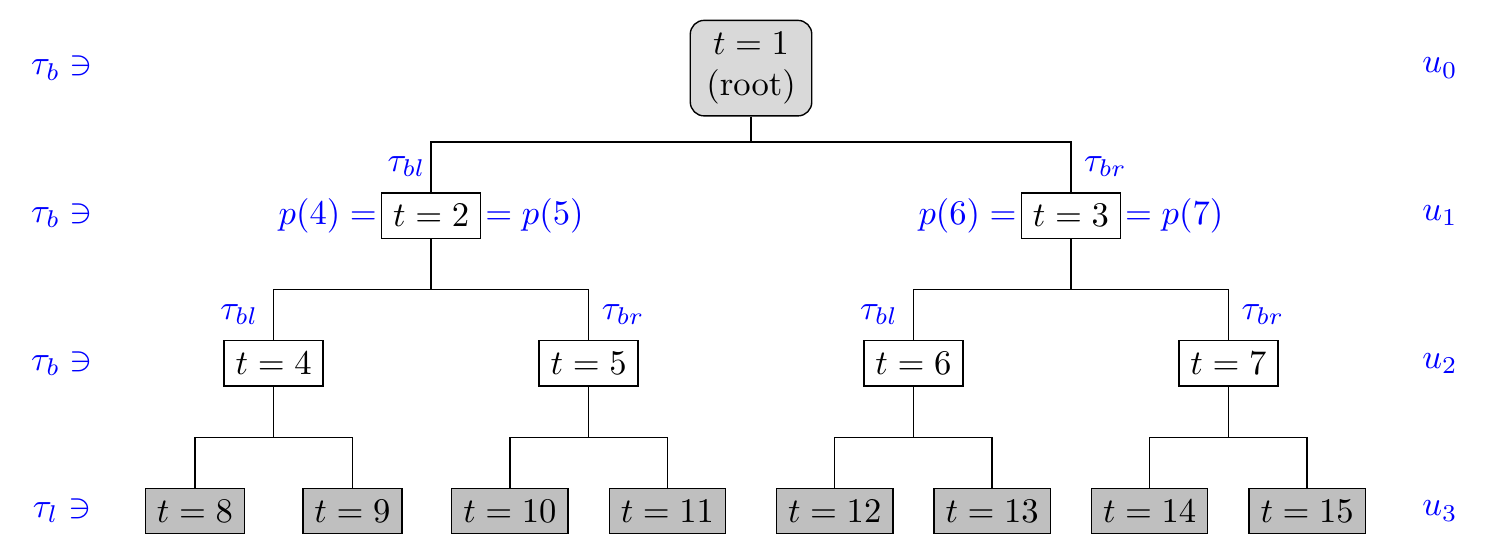}
\caption{Elements in a depth $D=3$ tree.\label{fig:tree}}
\end{figure}

In Figure \ref{fig:tree} we show the above mentioned elements in a $3$-depth tree.

Apart from the information about the topological structure of the tree, we also consider three regularization parameters that have to be calibrated in the validation process that allow us to find a trade-off between the different goals that we combine in our model: margin violation and classification errors of the separating splitting hyperplanes, correct classification at the leaf nodes and complexity of the tree. These parameters are the following:
 \begin{itemize}
\item[]  $c_1$:  unit misclassification cost at the leaf nodes.\\
 \item[] $c_2$: unit distance based misclassification errors for SVM splits.\\
\item[] $c_3$: unit cost for each splitting hyperplane introduced in the tree.
 \end{itemize}
 
The complete list of index sets and parameters used in our model are summarized in Table \ref{t:parameters}.
 
\begin{table}[h]
\begin{tabular}{rp{9cm}}\hline
$N = \{1, \ldots, n\}$ & Index set for the observations in the training sample.\\
$D$ & Maximal depth of the tree.\\
$T=2^{D+1}$ & maximal number of nodes in a $D$-depth tree.\\
$\tau = \{1, \ldots, T\}$ & Index set for the set of nodes of the tree.\\
$p(t)$ & Parent of node $t$, for $t\in \tau\backslash\{1\}$.\\
$\tau_b \in \tau$ & Branching nodes of the tree.\\
$\tau_l$ & Leaf nodes of the tree.\\
$\tau_{bl} \in \tau_b$ & Nodes that follow the  left branch on the path from their parent nodes. \\
$\tau_{br} \in \tau_b$ & Nodes whose right branch has been followed on the path from their parent nodes.\\
$u_s$ & Nodes at level $s$ of the tree, for $s=0, \ldots, D$.\\
$U = \{u_0, \ldots, u_d\}$ & Sets of levels of the tree.\\
 $c_1$ & Unit misclassification cost.\\
 $c_2$ & Unit distance based missclassification errors for SVM splits.\\
 $c_3$ & Unit cost for splitting hyperplanes.\\\hline
 \end{tabular}
 \caption{Index sets and parameters used in our model.\label{t:parameters}}
 \end{table}

\subsection{Variables}

Our model uses a set of decision and auxiliary variables that are described in Table \ref{table:vars}. We use both binary and continuous decision variables to model the MOCTSVM. The binary variables allow us to decide the allocation of observations to the decision tree nodes, or to decide whether a node is splited or not in the tree. The continuous variables allow us to determine the coefficients of the splitting hyperplanes or the misclassification errors (both in th SVM separations or at the leaf nodes). We also use auxiliary binary and integer variables that are useful to model adequately the problem.

\begin{table}[h]
\begin{tabular}{rp{9cm}}\hline
\multicolumn{2}{c}{Continuous Decision Variables}\\\hline
$\omega_t \in \mathbb{R}^p$ & Coefficients of the separating hyperplane of node $t$.\\
$\omega_{t_0} \in \R$ &Intercept of the separating hyperplane of node $t$.\\
$e_{it} \in \R_+$ &  Misclassification error of observation $i$ at  node $t$.\\
$\delta \in \R_+$ & Inverse of the minimum margin between splitting hyperplanes.\\\hline
\multicolumn{2}{c}{Binary Decision Variables}\\\hline
 $z_{it}\in\{0,1\}$ & Is one if observation $i$ belongs to node $t$ and zero otherwise.\\
$d_t\in\{0,1\}$ & Is one if a split is applied at node $t$ and zero otherwise.\\
\hline
\multicolumn{2}{c}{Auxiliary Variables}\\\hline
$L_t \in \Z_+$ & Number of misclassified observations at leaf node $t$.\\
$\alpha_{it} \in\{0,1\}$ & Is one if observation $i$ belongs to the reference fictitious class in node $t$ and zero otherwise.\\
$h_{it} \in \{0,1\}$ & Is one if observation $i$ is in node $t$ and lies on the positive half space of the hyperplane of node $t$, and zero otherwise.\\
$v_t \in \{0,1\}$ & Is one if not all observations in node $t$ lie on the positive half space of the hyperplane in node $t$ and zero otherwise.\\
$q_{kt} \in\{0,1\}$ & Is one if class $k$ is the most represented one in leaf node $t$ and zero otherwise.\\
\hline
 \end{tabular}
 \caption{Summary of the variables used in our model.\label{table:vars}}
 \end{table}
 
\begin{figure}[h]
 \includegraphics{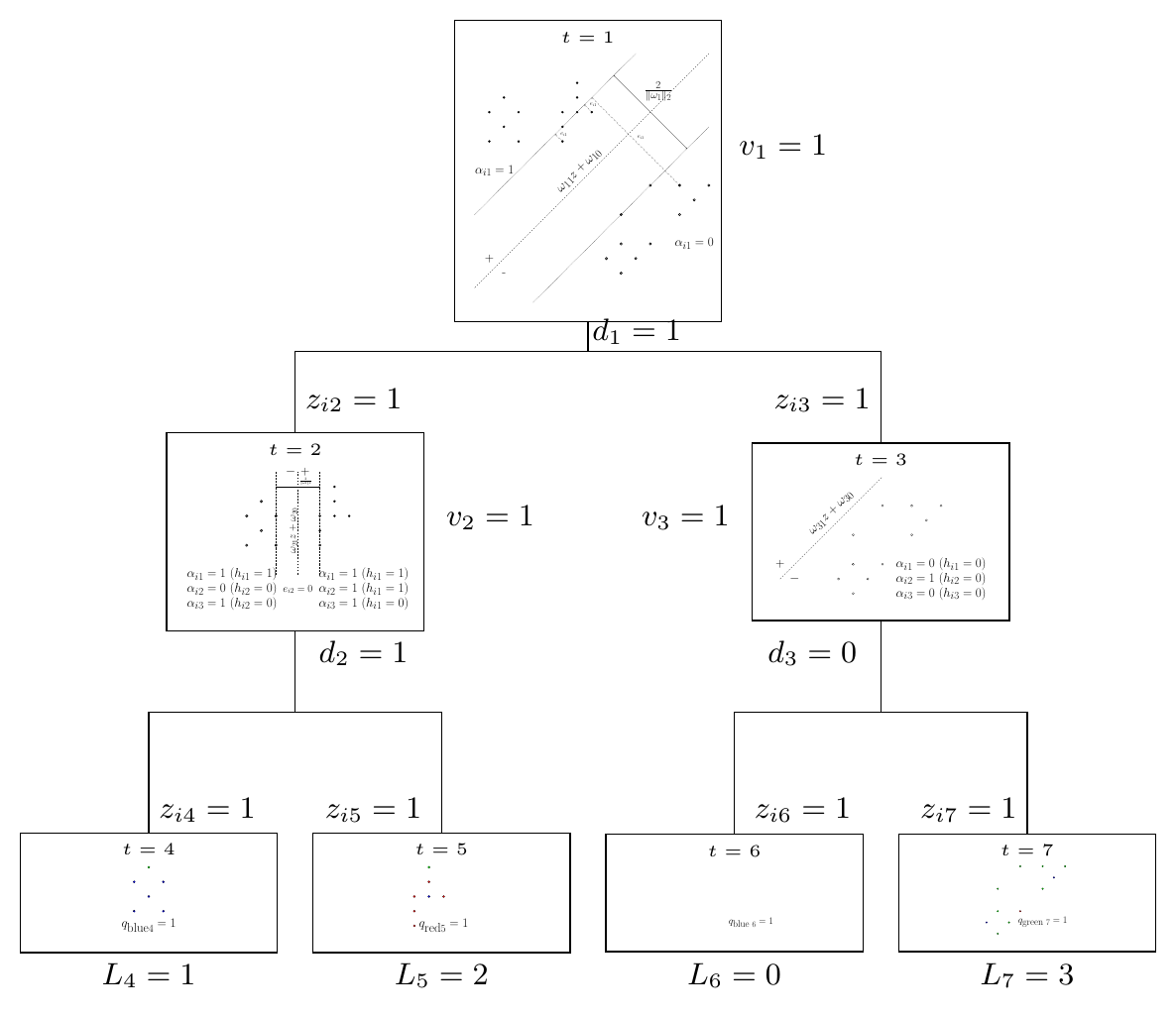}
 \caption{Illustration of the sets of variables used in our model in a toy example.\label{fig:tree2}}
 \end{figure}

In Figure \ref{fig:tree2} we illustrate the use of these variables in a feasible solution of a toy instance with three classes (red,blue and green).
 
The whole set of training observation is considered at the root node (node $t=1$). There, the original labels are ignored and to determine the \textit{fictitious} class of each observation a SVM-based hyperplane is constructed. Such a hyperplane is defined by the coefficients $\omega_1 \in \R^p$ and $\omega_{10}\in \R$ (hyperplane/line drawn with a dotted line in the picture) and it induces a margin separation ($\frac{2}{\|\omega_1\|_2}$) and misclassification errors $e_{i1}$. In the feasible solution drawn in the figure, only three observations induce positive errors (those that are classified either in the margin area or in the opposite side of the hyperplane). Such a hyperplane also determines the splitting rule for the definition of the children of that node. Since the node is  split ($d_1=1$), the observations that belongs to the \textit{positive side} of the hyperplane are assigned to the left node (node $t=2$) while those in the negative side are assigned to the right node (node $t=3$) through the $z$-variables. At node $t=2$, the same scheme is applied, that is, the hyperplane defined by $\omega_2$ is constructed, inducing SVM-based margin and errors and since $d_2=1$, also the splitting rule applies to define nodes $t=4$ and $t=5$. At node $t=2$,  one must control the observations in that node to quantify the  misclassifying errors, $e_{i2}$, only for those observations in the objective function. Specifically, we only account for these errors for the observations that belong to the node ($z_{i2}=1$) and either belong to the positive ($\alpha_{i2}=1$) or the negative ($\alpha_{i2}=0$) side of the hyperplane. Also, in order to control the complexity of the tree, the $h$-variables are used to know whether an observation belongs to the node and to the positive side of the SVM-hyperplane. If all observations in a node belong to the positive side of the hyperplane, the variable $v$ assumes the value $0$. Otherwise, in case $v$ takes value $1$, two situations are possible:  {1) there are observations in both sides of the hyperplane (as in node $t=2$) inducing a new split ($d_2=1$), and 2) all observations belong to the negative side (as in node $t=3$) determining that the tree is pruned at that node ($d_3=0$). 

Concerning the leaf nodes, node $t=2$ is split into nodes $t=4$ and $t=5$ and node $t=3$, which was decided to be no longer split, is fictitiously split in two leaf nodes, although one of them is empty and the other one receives all the observations of the parent node (node $t=3$). The allocation of any leaf node $\tau_l$ to a class is done through the $q$-variables (to the most popular class in the node or arbitrarily in case the node has no observations) and the number of misclassified observations is accounted for by the $L$-variables.}

\subsection{Objective Function}

As already mentioned, our method aims to construct classification trees with small misclassification errors at the leaf nodes, but at the same time with maximal separation between the classes with the SVM-based hyperplanes and minimum distance based errors. 

Using the variables described in the previous section, the four terms that are included in the objective functions are the following:

\begin{description}
\item[\bf Margins of the splitting hyperplanes:] The separating hyperplane of branching node $t \in \tau_b$ has margin $\frac{2}{\|\omega_t\|_2}$. Thus, our method aims to maximize the minimum of these margins. This is equivalent to minimize the maximum among the inverse margins $\{\frac{1}{2} \|\omega_t\|_2^2: t \in \tau_b\}$ which is represented by the auxiliary variable $\delta$.
\item[\bf Misclassification Errors at the leaf nodes:] Variable $L_t$ accounts for the number of misclassified observations in leaf node $t$, i.e., the number of observations that do not belong to the most represented class in that leaf node. These variables allow us to count the overall number of misclassified observations in the training sample. Therefore, the amount  to be minimized by the model is given by the following sum:
$$
c_1\dsum_{t\in \tau_l} L_t
$$
\item[\bf Distance-based Errors at branching nodes:] Each time a split is added to the tree, a SVM-based hyperplane in which the labels are assigned based on the global convenience of for the overall tree is incorporated. Thus, we measure, at each branching node in $\tau_b$, the distance-based errors incurred by the SVM classifier at that split. This amount is measured by the $e_{it}$ variables  and is incorporated to the model through the sum:
$$
c_2 \dsum_{i\in N} \dsum_{t\in \tau_b} e_{it}
$$
\item[\bf Complexity of the tree] The simplicity of the resulting tree is measured by the number of splits that are done in its construction. Since the $d_t$ variable tells us whether node $t$ is split or not, this term is accounted for in our model as:
$$
c_3 \dsum_{t\in \tau_b} d_t
$$
\end{description}

Summarizing, the overall objective  function of our model is:
\begin{equation}\label{eq:obj}\tag{${\rm OBJ}$}
\min\  \delta +  c_1\sum_{t\in\tau_l} L_t + c_2\sum_{i\in N}\sum_{t\in\tau_b}e_{it} + c_3\sum_{t\in\tau_b} d_t.
\end{equation}

Note that the coefficients $c_1$, $c_2$ and $c_3$ trade-off the misclassification of the training sample, the separation between classes and the complexity of the tree, respectively. These parameters should be carefully calibrated in order to construct \textit{simple} decision trees with high predictive power, as can be seen in our computational experiments.

\subsection{Constraints}

The requirements on the relationships between the variables and the rationale of our model are described through the following constraints that define the mathematical programming model.

First of all, in order to adequately represent the maximum among the inverse margins of the sppliting hyperplanes, we require:
\begin{equation}\label{c:1}\tag{${\rm C1}$}
\delta \geq \frac{1}{2} \|\omega_t\|_2^2, \forall t \in \tau_b.
\end{equation}
Next, we impose how the splits are performed in the tree. To this end, we need to know which observations belong to a certain node $t$ ($z$-variable) and how these observations are distributed with respect to the two fictitious classes to be separated ($\alpha$-variables). Gathering all these elements together, we use the following constraints to define the splits of the decision tree:
\begin{align}
& \omega_t'x_i + \omega_{t0} \geq 1 - e_{it} -(1-h_{it}) &\forall i\in N, t\in\tau_b,\label{c:2a}\tag{${\rm C2a}$}\\
& \omega_t'x_i + \omega_{t0} \leq -1 + e_{it} + (1-z_{it}+\alpha_{it}) &\forall i\in N,t\in\tau_b.\tag{${\rm C2b}$}\label{c:2b}
\end{align}
According to this, constraint \eqref{c:2a} is activated just in case the  observation $i$ belongs to the reference class and it is in node $t$ ($h_{it}=1$). On the other hand, \eqref{c:2b} is activated if $i$ is allocated to node $t$ ($z_{it}=1$) but it does not belong to the reference class ($\alpha_{it}=0$). Therefore, the reference class is located on the positive half space of hyperplane $\mathcal{H}_t$, while the other class is positioned in the negative half space, and at the same time, margin violations are regulated by the $e_{it}$ variables.

To ensure the correct behaviour of the above constraints, we must correctly define the $z_{it}$ variables. First, it is required that  each observation belongs to exactly one node per level in the tree. This can be easily done by adding the usual assingment constraints to the problem  at each of the levels, $u\in U$, of the tree:
\begin{align}\label{c:3}
& \sum_{t\in u}z_{it} = 1 &\forall i\in N, u\in U. \tag{${\rm C3}$}
\end{align}
Furthermore, we should enforce that if observation $i$ is in node $t$ ($z_{it}=1$), then observation $i$ must also be in the parent node of $t$, $p(t)$ ($z_{ip(t)}=1$), and also observation $i$ can not be in node $t$ if it is not in its parent node ($z_{ip(t)}=0 \Rightarrow z_{it}=0$). These implications can be obtained by means of the following constraints:
\begin{align}\label{c:4}
& z_{it} \leq z_{ip(t)} &\forall i\in N, t=2,\ldots, T.\tag{${\rm C4}$}
\end{align}
Nevertheless, the way observations descend through the tree needs a further analysis, since at this point they could just randomly define a path in the tree. Whenever an observation $i$ is in the positive half space of the splitting hyperplane at  node $t$, $\mathcal{H}_t$, this observation should follow the right branch connecting to the child node of $t$. Otherwise, in case $i$ is on the negative half space, it should follow the left branch. The knowledge on the side of the splitting hyperplane where an observation belongs to is encoded in the $\alpha$-variables. Then, in case $i$ lies on the positive half space of $\mathcal{H}_t$, $\alpha_{it}$ will never be equal to zero since it would lead to a value of $e_{it}$ greater than one, while $e_{it}<1$ is guaranteed in case $\alpha_{it}=1$.

With the above observations, the constraints that assure the correct construction of the splitting hyperplanes with respect to the side of them where the observations belong to are the following:
\begin{align}
& z_{ip(t)} - z_{it} \leq \alpha_{ip(t)} & \forall i\in N, t \in \tau_{bl}, \label{c:5a}\tag{${\rm C5a}$}\\
& z_{ip(t)} - z_{it} \leq 1-\alpha_{ip(t)} & \forall i\in N, t \in \tau_{br} \label{c:5b}.\tag{${\rm C5b}$}
\end{align}
Constraints \eqref{c:5a} assure that if observation $i$ is on the parent node of an even node $t$ ($z_{ip(t)}=1$), and $i$ lies on the negative half space of $\mathcal{H}_{p(t)}$ ($\alpha_{ip(t)}=0$), then $z_{it}$ is enforced to be equal to one. As a result, $\alpha_{ip(t)}=0$ forces observation $i$ to take the left branch in node $t$. Note that in case $z_{ip(t)}=1$,  and at the same time observation $i$ is not in the left child node of $t$ ($z_{it}=0$ for $i\in\tau_{bl}$), then  $\alpha_{ip(t)}=1$,  which means that observation $i$ lies on the positive half space of $\mathcal{H}_{p(t)}
$.  Constraints \eqref{c:5b} are analogous to \eqref{c:5a} but allowing to adequately represent right branching nodes. 

Moreover, two additional important elements need to be incorporated to complete our model:  the tree complexity and the correct definition of misclassified observations. Note that in usual Optimal Classification Trees that do not use SVM-based splits, the complexity can be easily regulated by just imposing $\|\omega_t\|_2^2 \leq M d_t$ (for a big enough $M$ constant) in all the branch nodes, since in case a node is no further branched ($d_t=0$), the coefficients of the splitting hyperplane are set to zero. However, in our case, in which the splitting hyperplanes are SVM-based hyperplanes, these constraints are in conflict with constraints \eqref{c:2a} and \eqref{c:2b}, since in case $d_t = 0$ (and therefore $\omega_t=0$) it would not only imply that the coefficients $\omega_t$ are equal to zero, but also that the distance based errors would be set to the maximum value of $1$, i.e., $e_{it}=1$ for every observation $i$ in the node, even though these errors would not make any sense since observations would not be separated at the node. To overcome this issue, we consider the auxiliary binary variables $h_{it}=z_{it} \alpha_{it}$ ($h_{it}$ takes value 1 if observation $i$ belongs to node $t$ and lies in the positive half-space of the splitting hyperplane applied at node $t$) and $v_t$ (that takes value zero in case all the points in the node belong to the positive halfspace and one otherwise). The variables are adequatelly defined if the following constraints are incorporated to the model:
\begin{align}
& h_{it} \geq z_{it} + \alpha_{it} -1,& \forall i\in N, t\in\tau_b, \label{c:6a}\tag{${\rm C6a}$}\\
& h_{it} \leq z_{it} - \alpha_{it} +1,& \forall i\in N, t\in\tau_b, \label{c:6b}\tag{${\rm C6b}$}\\
&\sum_{i\in N} (z_{it}-h_{it}) \leq n v_t, &\forall t\in\tau_b, \label{c:6c}\tag{${\rm C6c}$}\\
&\sum_{i\in N} h_{it} \leq n (1+d_t-v_t),  &\forall  t\in\tau_b , \label{c:6d}\tag{${\rm C6d}$}
\end{align}
where constraints \eqref{c:6a} and \eqref{c:6b} are the linearization of the bilinear constraint $h_{it}=z_{it} \alpha_{it}$. On the other hand, Constraints \eqref{c:6c} assure that in case $v_t=0$, then all observations in node $t$ belong to the positive halfspace of $\mathcal{H}_t$, and constraints  \eqref{c:6d} assure  that if $v_t=1$ and the tree is pruned at node $t$ ($d_t=0$), then those observations allocated to node $t$ are placed in the negative halfspace defined by the splitting hyperplane. Thus, it implies that $d_t$ takes value one if and only if the observations in node $t$ are separated by $\mathcal{H}_t$, and therefore producing an effective split at the node.

Finally, in order to adequately represent the $L_t$ variables (the ones that measure the number of misclassified observations at the leaf nodes) we use the constraints already incorporated in the OCT-H model in \cite{bertsimas2017optimal}. On the one hand, we assign each leaf node to a single class (the most popular class of the observations that belong to that node). We use the binary variable $q_{kt}$ to check whether leaf node $t\in \tau_l$ is assigned to class $k=1, \ldots, K$. The usual assignment constraints are considered to assure that each node is assigned to exactly one class:
\begin{align}
& \sum_{k=1}^K q_{kt} = 1,  &\forall  t \in\tau_l.  \label{c:7}\tag{${\rm C7}$}
\end{align}

The correct definition of the variable $L_t$ is then guaranteed by the following set of constraints:
\begin{align}
& L_t \geq  \sum_{i\in N} z_{it}  - \sum_{i\in N} Y_{ik}z_{it} -  n (1-q_{kt}), &\forall  k=1,\ldots, K, t \in\tau_l, \label{c:8} \tag{${\rm C8}$}
\end{align}
These constraints are activated if and only if $q_{kt}=1$, i.e., if observations in node $t$ are assigned to class $k$. In such a case, since $L_t$ is being minimized in the objective function, $L_t$ will be determined by the number of training observations in node $t$ except those whose label is $k$, i.e., the number of missclasified observations in node $t$ according to the $k$-class assignment.

Observe that the constant $n$ in \eqref{c:8} can be decreased and fixed to the maximum number of misclassified observations in the training sample. This number coincide with the difference between the number of observations in the training sample ($n$) and  the number of observations in the most represented class in the sample.

Summarizing the above paragraphs, the MOCTSVM can be formulated as the following MINLP problem:
\begin{align} 
\min &\;\;\; \delta +  c_1\sum_{t\in\tau_l }L_t + c_2\sum_{i\in N}\sum_{t\in \tau} e_{it}+  c_3\sum_{t\in \tau} d_t \label{eq:obj}\tag{${\rm OBJ}$}\\
\mbox{s.t.} & \;\;\label{c:1}\tag{${\rm C1}$} \delta \geq \frac{1}{2} \|\omega_t\|, &\forall t \in \tau_b,\\
& \omega_t'x_i + \omega_{t0} \geq 1 - e_{it} -(2-z_{it}-\alpha_{it}), &\forall i\in N, t\in\tau_b,\label{c:2a}\tag{${\rm C2a}$}\\
& \omega_t'x_i + \omega_{t0} \leq -1 + e_{it} + (1-z_{it}+\alpha_{it}), &\forall i\in N,t\in\tau_b,\tag{${\rm C2b}$}\label{c:2b}\\
& \sum_{t\in u}z_{it} = 1, &\forall i\in N, u\in U, \tag{${\rm C3}$}\label{c:3}\\
& z_{it} \leq z_{ip(t)}, &\forall i\in N, t=2,\ldots, T,\label{c:4a}\tag{${\rm C4}$}\\
& z_{ip(t)} - z_{it} \leq \alpha_{ip(t)}, & \forall i\in N, t \in \tau_{bl}, \label{c:5a}\tag{${\rm C5a}$}\\
& z_{ip(t)} - z_{it} \leq 1-\alpha_{ip(t)}, & \forall i\in N, t \in \tau_{br}, \label{c:5b}.\tag{${\rm C5b}$}\\
& h_{it} \geq z_{it} + \alpha_{it} -1,& \forall i\in N, t\in\tau_b, \label{c:6a}\tag{${\rm C6a}$}\\
& h_{it} \leq z_{it} - \alpha_{it} +1,& \forall i\in N, t\in\tau_b, \label{c:6b}\tag{${\rm C6b}$}\\
&\sum_{i\in N} (z_{it}-h_{it}) \leq n v_t, &\forall t\in\tau_b, \label{c:6c}\tag{${\rm C6c}$}\\
&\sum_{i\in N} h_{it} \leq n (1+d_t-v_t),  &\forall  t\in\tau_b, \label{c:6d}\tag{${\rm C6d}$}\\
& \sum_{k=1}^K q_{kt} = 1,  &\forall  t \in\tau_l,  \label{c:7}\tag{${\rm C7}$}\\
& L_t \geq  \sum_{i\in N} z_{it}  - \sum_{i\in N} Y_{ik}z_{it} - n (1-q_{kt}), &\forall  k=1,\ldots, K, t \in\tau_l, \label{c:8} \tag{${\rm C8}$}\\
& e_{it} \in\mathbb{R}^+,  \alpha_{it}, h_{it} \in \{0,1\}, & \forall i\in N, t\in\tau_b,\nonumber\\
& z_{it} \in \{0,1\}, & \forall i\in N,t=1,\ldots, T, \nonumber\\
& q_{kt}\in \{0,1\}, & \forall k=1\ldots,K,t\in\tau_l,\nonumber\\
&\omega_t\in\mathbb{R}^p, \omega_{t0} \in\mathbb{R}, d_t \in\{0,1\}, & \forall t=1,\ldots,T.\nonumber
\end{align}

\subsection{Strengthening the model}

The MINLP formulation presented above is valid for our MOCTSVM model. However, it is a computationally costly problem, and although it can be solved by most of the off-the-shelf optimization solvers (as Gurobi, CPLEX or XPRESS), it is able to solve optimally only small to medium size instances. To improve its performance, the problem can be strengthen by means of valid inequalities which allows one to reduce the gap between the continuous relaxation of the problem and its optimal integer solution, being then able to solve larger instances in smaller CPU times. In what follows we describe some of these inequalities that we have incorporated to the MINLP formulation:
\begin{itemize}
\item If observations $i$ and $i'$ belongs to different nodes, they cannot be assigned to the same node for the remainder levels of the tree:
$$
z_{is} + z_{i's} \leq z_{it} + z_{i't}, \forall t \in u, s \in u' \text{$u \leq u'$}
$$
\item If leaf nodes $t$ and $s$ are the result of proper splitting hyperplanes, then, both nodes cannot be assigned to the same class:
$$
q_{kt} + q_{ks} \leq 2-d_{p(t)}, \forall t, s  =2, \ldots, T (t\neq s) \text{ with } p(t)=p(s), k=1, \ldots, K.
$$
\item Variable $\alpha_{it}$ is enforced to take value $0$ in case $z_{it}=0$:
$$
\alpha_{it} \leq z_{it}, \forall i \in N, t \in \tau_b.
$$
\item Variable $h_{it}$ is not allowed to take value one if  $\alpha_{it}$ takes value zero:
\begin{align*}
h_{it} \leq \alpha_{it}, \forall i \in N, t \in \tau_b.
\end{align*}
\item There should be at least a leaf node to which each class is assigned to (assuming that each class is represented in the training sample). It also implies that the number of nodes to which a class is assigned is bounded as:
\begin{align*}
1 \leq \dsum_{t\in \tau_l} q_{kt} \leq 2^D - 1, \forall k=1, \ldots, K.
\end{align*}
\end{itemize}

In order to reduce the dimensionality and also to avoid symmetries of the MINLP problem, one can also apply some heuristic strategies to fix the values of some of the binary variables in a preprocessing phase. For instance, we choose $i_0 \in \arg\displaystyle\max_{i\in N} \mid\{i' \in N: \|a_{i}-a_{i'}\|\leq \varepsilon \text{ and } y_i=y_{i'}\}\mid$ , that is, the observation with a maximum amount of observations in the same class \textit{close enough to} it. Then, we fix to one all the variables $z_{it_0}$ with $i \in \{i' \in N: \|a_{i}-a_{i'}\|\leq \varepsilon \text{ and } y_i=y_{i'}\}$, being $t_0$ the first left leaf node of the tree (and fixing to zero the allocation of these points to the rest of the leaf nodes). Analogously, we fix also to one all the $z$-variables allocating observations to the last right leaf node of the tree for a subset of observations in the same class which are far enough from $i_0$, i.e., $z_{it_{f}}=1$ for all $i \in \{i' \in N: \|a_{i_f}-a_{i'}\|\leq \varepsilon \text{ and } y_{i_f}=y_{i'}\}$ where $i_f = \arg\displaystyle\max_{i\in N} \|a_{i_0}-a_i\|$  (and fixing to zero the allocation of these points to the rest of the leaf nodes).

\section{Experiments}\label{sec:5}

In order to analyze the performance of this new methodology we have run a series of experiments among different real datasets from UCI Machine learning Repository \cite{uci}. We have chosen twelve datasets with number of classes between two and seven. The dimension of these problems is reported in Table \ref{table:1} by the tuple $(n: \text{number of observations}, p:\text{number of features}, K: \text{number of classes})$.

We have compared the MOCTSVM model with three other Classification Tree-based methodologies, namely CART, OCT and OCT-H. The maximum tree depth, $D$, for all the models was equal to $3$, and the minimum number of observations per node in CART, OCT and OCT-H was equal to the $5\%$ of the training size.

We have performed, for each instance a 5-fold cross validation scheme, i.e., datasets have been splited into five random train-test partitions where one of the folds is used to build the model and the remaining are used to measure the accuracy of the predictions. Moreover, in order to avoid taking advantage of beneficial initial partitions, we have repeated the cross-validation scheme five times for all the datasets. 

The CART method was coded in \texttt{R} using the \texttt{rpart} library. On the other hand, MOCTSVM, OCT and OCT-H were coded in \texttt{Python} and solved using the optimization solver \texttt{Gurobi} 8.1.1. All the experiments were run on a PC Intel Xeon E-2146G processor at 3.50GHz and 64GB of RAM. A time limit of 300 seconds was set for training the training folds.  Although not all the problems were optimally solved within the time limit, as can be observed in Table \ref{table:1}, the results obtained with our model already outperform the other methods.

In order to calibrate the parameters of the different models that regulate the complexity of the tree, we have used different approaches. On the one hand, for CART and OCT, since the maximum number of nodes for such a depth is $2^D-1 = 7$,  one can search for the tree with best complexity by searching in the grid $\left\lbrace 1,\ldots, 2^D-1\right\rbrace$ of possible active nodes.  For OCT-H, we search the complexity regularization factor in the grid ${\left\lbrace 10^i: i=-5,\ldots , 5  \right\rbrace}$. Finally, in MOCTSVM we used the same grid  ${\left\lbrace 10^i: i=-5,\ldots , 5  \right\rbrace}$ for $c_1$ and $c_2$, and ${\left\lbrace 10^i: i=-2,\ldots , 2  \right\rbrace}$ for $c_3$.

In Table \ref{table:1} we report the results obtained in our experiments for all the models. The first column of the table indicates the identification of the dataset (together with its dimensionality). Second, for each of the methods that we have tested, we report the obtained average test accuracy and the standard deviation. We have highlighted in bold the best average test accuracies obtained for each dataset.

As can be observed, our method clearly outperforms in most of the instances the rest of the methods in terms of accuracy. Clearly, our model is designed  to construct Optimal Classification Trees with larger separations between the classes, which results in better accuracies in the test sample. The datasets \textit{Australian} and \textit{BalanceScale} obtain their better results with OCT-H, but, as can be observed, the differences with respect the rest of the methods are tiny (it is the result of correctly classifying in the test sample just a few more observations  than the rest of the methods). In that case, our method gets an accuracy almost as good as OCT-H. In the rest of the datasets, our method consistently gets better classifiers and for instance for \textit{Dermatology} the difference with respect to the best classifiers among the others ranges in $[4\%, 19\%]$, for  \textit{Parkinson} the accuracy with our model is at least $6\%$ better than the rest, for Wine we get $5\%$ more accuracy than OCTH and $10\%$ more than CART and for  Zoo the accuracy of our model is more than $17\%$ greater than the one obtained with CART. 

Concerning the variability of our method, the standard deviations reported in Table \ref{table:1} show  that our results are, in average, more stable than the others, with small deviations with respect to the average accuracies. This behaviour differs from the one observed in CART or OCT, where larger deviations are obtained, implying that the accuracies highly depends of the test folder where the method is applied.

\begin{table}[h!]
\hspace*{-1cm}\begin{tabular}{c|ccccc}
\multicolumn{1}{c|}{} & CART & OCT & OCT-H & MOCTSVM & Diff\\
\hline 
\begin{tabular}[c]{@{}c@{}}Australian \\(690,14,2)\end{tabular} & 85.54 $\pm$ 0.81  & 85.22 $\pm$ 1.27 & {\bf 85.65 $\pm$ 1.02} & 85.27 $\pm$ 1.11 & -0.38 $\pm$ 0.63 \\
\begin{tabular}[c]{@{}c@{}}BalanceScale \\(625,4,3)\end{tabular} & 69.55 $\pm$ 1.76  & 73.30 $\pm$ 1.20 & {\bf 90.43 $\pm$ 1.07 } & 89.53 $\pm$ 1.28 &  -0.90 $\pm$ 1.69\\
\begin{tabular}[c]{@{}c@{}}Banknote \\(1372,5,2)\end{tabular} &    89.27 $\pm$ 0.95 & 88.50 $\pm$ 1.17 & 98.89 $\pm$ 0.33 &  {\bf 98.91 $\pm$ 0.46 } & 0.02 $\pm$ 0.43\\
\begin{tabular}[c]{@{}c@{}}BreastCancer \\(683,9,2)\end{tabular} & 92.69 $\pm$ 1.01  & 94.16 $\pm$ 0.54 & 95.10 $\pm$ 1.26 & {\bf 96.27 $\pm$ 0.64} & 1.17 $\pm$ 1.39     \tabularnewline
\begin{tabular}[c]{@{}c@{}}Dermatology \\(358,34,6)\end{tabular} & 75.69 $\pm$ 3.60  & 77.82 $\pm$ 4.34 & 91.41 $\pm$ 2.83 & {\bf 95.39 $\pm$ 1.47 } & 3.98 $\pm$ 2.74    \tabularnewline
\begin{tabular}[c]{@{}c@{}}Heart \\(294,13,5)\end{tabular} &  64.37 $\pm$ 1.48  & 65.14 $\pm$ 1.57  & 64.30 $\pm$ 1.79 &  {\bf 66.41 $\pm$ 1.54}   & 1.26 $\pm$ 1.38   \tabularnewline
\begin{tabular}[c]{@{}c@{}}Iris \\(150,4,3)\end{tabular} &  94.26 $\pm$ 1.90 & 95.37 $\pm$ 0.97 & 95.64 $\pm$ 1.46  & {\bf 95.72 $\pm$ 1.79}  & 0.08 $\pm$ 1.70    \tabularnewline
\begin{tabular}[c]{@{}c@{}}Parkinson \\(240,40,2)\end{tabular} & 72.29 $\pm$ 4.05 & 73.53 $\pm$ 2.26 & 74.92 $\pm$ 3.01 & {\bf 80.83 $\pm$ 1.89}  & 5.91 $\pm$ 3.06    \tabularnewline
\begin{tabular}[c]{@{}c@{}}Seeds \\(210,7,3)\end{tabular} &  86.36 $\pm$ 4.02 & 88.52 $\pm$ 2.69 & 91.12 $\pm$ 2.99 &{\bf 92.98 $\pm$ 1.82} & 1.85 $\pm$ 2.23 \tabularnewline
\begin{tabular}[c]{@{}c@{}}Teaching \\(150,5,3)\end{tabular} & 41.91 $\pm$ 5.64 & 48.35 $\pm$ 3.80  & 48.09 $\pm$ 2.92 & {\bf 48.62 $\pm$ 3.37} & 0.26 $\pm$ 4.60    \tabularnewline
\begin{tabular}[c]{@{}c@{}}Thyroid \\(215,5,3)\end{tabular} & 89.77 $\pm$ 2.37 & 92.43 $\pm$ 2.12 & 92.46 $\pm$ 2.49 & {\bf 94.57 $\pm$ 2.08}  & 2.11 $\pm$  3.10    \tabularnewline
\begin{tabular}[c]{@{}c@{}}Wine \\(178,13,3)\end{tabular}    & 84.52 $\pm$ 2.66 & 92.22 $\pm$ 3.41 &  89.35 $\pm$ 3.71 & {\bf 94.13 $\pm$ 1.78} & 1.90 $\pm$ 3.37 \tabularnewline
\begin{tabular}[c]{@{}c@{}}Zoo\\(101,16,7)\end{tabular} & 74.96 $\pm$ 5.79 & 87.75 $\pm$ 1.99 & 89.11 $\pm$ 2.58 & {\bf 92.31 $\pm$ 2.15} & 3.20 $\pm$ 2.95     \tabularnewline
\hline 
\end{tabular}
\caption{Average accuracies ($\pm$ standard deviations) obtained in our computational experiments.\label{table:1}}
\end{table}
\section{Conclusions and Further Research}\label{sec:6}
We have presented in this paper a novel methodology to construct classifiers for multiclass instances by means of a Mathematical Programming model. The proposed method outputs a classification tree were the splits are based on SVM-based hyperplanes. At each branch node of the tree, a binary SVM hyperplane is constructed in which the observations are classified in two fictitious classes (the original classes are ignored in all the splitting nodes), but the global goodness of the tree is measured at the leaf nodes, where misclassification errors are minimized. Also, the model minimizes the complexity of the tree together with the two elements that appear in SVM-approaches: margin separation and distance-based misclassifying errors.  We have run an extensive battery of computational experiments that shows that our method outperforms most of the Decision Tree-based methodologies both in accuracy and stability. 

Future research lines on this topic include the analysis of nonlinear splits when branching  in  MOCTSVM, both using kernel tools derived from SVM classifiers or  specific families of nonlinear separators. This approach will result into more flexible classifiers able to capture the nonlinear trends of many real-life datasets. Additionally, we also plan to incorporate  features selection in our method in order to construct  highly predictive but also more interpretable classification tools.

\section*{Acknowledgements}

This research has been partially supported by Spanish Ministerio de Ciencia e Innovación, Agencia Estatal de Investigación/FEDER grant number PID2020-114594GBC21, Junta de Andalucía projects P18-FR-1422 and projects FEDERUS-1256951, BFQM-322-UGR20, CEI-3-FQM331 and NetmeetData-Ayudas Fundación BBVA a equipos de investigación científica 2019. The first author was also partially supported by the IMAG-Maria de Maeztu grant CEX2020-001105-M /AEI /10.13039/501100011033.

\bibliographystyle{acm}
\bibliography{moctsvm_bib}

\end{document}